\documentclass[leqno]{llncs}

\usepackage{amsfonts}
\usepackage{amsmath}
\usepackage{amssymb}
\usepackage{comment}
\usepackage{expdlist}
\usepackage[dvips]{graphicx}
\usepackage{subfigure}

\newcommand{\set}[1]{\left \{ #1 \right \}}                     
\newcommand{\setst}[2]{\left \{ #1 \mid #2 \right \}}           
\newcommand{\abs}[1]{\left| #1 \right|}

\renewcommand{\phi}{\varphi}

\renewcommand{\bar}{\overline}
\renewcommand{\epsilon}{\varepsilon}

\newcommand{\calC}{\mathcal{C}}

\newcommand{\calP}{\mathcal{P}}

\newcommand{\calT}{\mathcal{T}}

\newcommand{\deltain}{\delta^{\rm in}}
\newcommand{\deltaout}{\delta^{\rm out}}

\renewcommand{\vec}[1]{\overrightarrow{#1}}
\newcommand{\val}[1]{\mathop{\rm val}\nolimits ( #1 )}

\newenvironment{numitem}{\refstepcounter{equation}\begin{enumerate}%
\item[(\theequation)]$\quad$}{\end{enumerate}}

\newcommand{\refeq}[1]{(\ref{eq:#1})}                
\newcommand{\reffig}[1]{Fig.~\ref{fig:#1}}           
\newcommand{\refth}[1]{Theorem~\ref{th:#1}}          
\newcommand{\reflm}[1]{Lemma~\ref{lm:#1}}            
\newcommand{\refsec}[1]{Section~\ref{sec:#1}}        
\newcommand{\refsubsec}[1]{Subsection~\ref{subsec:#1}}

\renewenvironment{proof}{\par\noindent%
{\bf Proof.\par\nopagebreak}}{\unskip\nobreak\enskip$\square$\par\bigskip}

\makeindex

\begin{document}

\title{A Linear Time Algorithm for \\ Finding Three Edge-Disjoint Paths \\ in Eulerian Networks}

\institute
{
    Moscow State University
}

\author
{
   Maxim A. Babenko
   \thanks{
   	Email: \texttt{max@adde.math.msu.su}.
        Supported by RFBR grant 09-01-00709-a.
   },
   Ignat I. Kolesnichenko \thanks{Email: \texttt{ignat1990@gmail.com}},
   Ilya P. Razenshteyn \thanks{Email: \texttt{ilyaraz@gmail.com}}
}

\maketitle

\begin{abstract}
    Consider an undirected graph $G = (VG, EG)$ and a set of six \emph{terminals}
    $T = \set{s_1, s_2, s_3, t_1, t_2, t_3} \subseteq VG$.
    The goal is to find a collection $\calP$ of
    three edge-disjoint paths $P_1$, $P_2$, and $P_3$,
    where $P_i$ connects nodes~$s_i$ and~$t_i$ ($i = 1, 2, 3$).

    Results obtained by Robertson and Seymour by graph minor techniques
    imply a polynomial time solvability of this problem.
    The time bound of their algorithm is $O(m^3)$
    (hereinafter we assume $n := \abs{VG}$, $m := \abs{EG}$, $n = O(m)$).

    In this paper we consider a special, \emph{Eulerian} case of $G$ and $T$.
    Namely, construct the \emph{demand graph} $H = (VG, \set{s_1t_1, s_2t_2, s_3t_3})$.
    The edges of $H$ correspond to the desired paths in $\calP$.
    In the Eulerian case the degrees of all nodes in the (multi-) graph $G + H$
    ($ = (VG, EG \cup EH)$) are even.

    Schrijver showed that, under the assumption of Eulerianess,
    cut conditions provide a criterion for the existence of $\calP$.
    This, in particular, implies that checking for existence of $\calP$
    can be done in $O(m)$ time.
    Our result is a combinatorial $O(m)$-time algorithm
    that constructs $\calP$ (if the latter exists).
\end{abstract}

\section{Introduction}
\label{sec:intro}

We shall use some standard graph-theoretic notation through the paper.
For an undirected graph $G$, we denote its sets of nodes and edges by $VG$
and $EG$, respectively. For a directed graph, we speak of arcs rather than
edges and denote the arc set of $G$ by $AG$. A similar notation
is used for paths, trees, and etc.
We allow parallel edges and arcs but not loops in graphs.

For an undirected graph $G$ and $U \subseteq VG$, we write $\delta_G(U)$ to
denote the set of edges with exactly one endpoint in $U$.
If $G$ is a digraph then the set of arcs entering (resp. leaving) $U$
is denoted by $\deltain_G(U)$ and $\deltaout_G(U)$.
For a graph $G$ and a subset $U \subseteq VG$, we write $G[U]$ to denote
the subgraph of $G$ induced by $U$.

Let $G$ be an undirected graph. Consider six nodes $s_1, s_2, s_3, t_1, t_2, t_3$ in $G$.
These nodes need not be distinct and will be called \emph{terminals}.
Our main problem is as follows:
\begin{numitem}
\label{eq:three_pairs}
    Find a collection of three edge-disjoint paths
    $P_1$, $P_2$, $P_3$, where $P_i$ goes from $s_i$ to $t_i$ (for $i = 1, 2, 3$).
\end{numitem}

Robertson and Seymour \cite{RS-95} developed sophisticated graph minor techniques
that, in particular, imply a polynomial time solvability of the above problem.
More specifically, they deal with the general case where
$k$ pairs of terminals $\set{s_1, t_1}, \ldots, \set{s_k, t_k}$ are given
and are requested to be connected by paths.
These paths are required to be \emph{node-disjoint}.
The edge-disjoint case, however, can be easily simulated by considering the line graph of $G$.
For fixed $k$, the running time of the algorithm of Robertson and Seymour
is cubic in the number of nodes (with a constant heavily depending on $k$).
Since after reducing the edge-disjoint case to the node-disjoint one
the number of nodes becomes $\Theta(m)$, one gets an algorithm of time complexity $O(m^3)$
(where, throughout the paper, $n := \abs{VG}$, $m := \abs{EG}$;
moreover, it is assumed that $n = O(m)$).
If $k$ is a part of input, then it was shown by Marx \cite{marx-04} that finding $k$ edge-disjoint paths is NP-complete
even in the Eulerian case.

We may also consider a general \emph{integer multiflow problem}.
To this aim, consider an arbitrary (multi-)graph $G$ and
also an arbitrary (multi-)graph $H$ obeying $VH = VG$.
The latter graph $H$ is called the \emph{demand graph}.
The task is to find a function~$f$ assigning each edge $uv \in EH$
a $u$--$v$ path $f(uv)$ in $G$ such that the resulting paths $\setst{f(e)}{e \in EH}$ are edge-disjoint.
Hence, the edges of $H$ correspond to the paths in the desired solution.
By a \emph{problem instance} we mean the pair $(G, H)$.
An instance is \emph{feasible} if the desired collection of
edge-disjoint paths exists; \emph{infeasible} otherwise.

In case of three terminal pairs one has $H = (VG, \set{s_1t_1, s_2t_2, s_3t_3})$.
We can simplify the problem and get better complexity results by introducing some
additional assumptions regarding the degrees of nodes in $G$.
Put $G + H := (VG, EG \cup EH)$.
Suppose the degrees of all nodes in $G + H$ are even;
the corresponding instances are called \emph{Eulerian}.

As observed by Schrijver, for the Eulerian case there exists
a simple feasibility criterion.
For a subset $U \subseteq VG$ let $d_G(U)$ (resp. $d_H(U)$)
denote $\abs{\delta_G(U)}$ (resp. $\abs{\delta_H(U)}$).

\begin{theorem}[\cite{sch-03}, Theorem~72.3]
\label{th:feasibility_criterion}
    An Eulerian instance $(G,H)$ with three pairs of terminals
    is feasible if and only if $d_G(U) \ge d_H(U)$ holds
    for each $U \subseteq VG$.
\end{theorem}

The inequalities figured in the above theorem are called \emph{cut conditions}.
In a general problem (where demand graph $H$ is arbitrary)
these inequalities provide necessary (but not always sufficient) feasibility requirements.

For the Eulerian case, the problem is essentially equivalent to constructing two paths (out
of three requested by the demand graph). Indeed, if edge-disjoint paths $P_1$ and $P_2$
(where, as earlier, $P_i$ connects $s_i$ and $t_i$, $i = 1, 2$) are found, the remaining
path $P_3$ always exists. Indeed, remove the edges of $P_1$ and $P_2$ from~$G$.
Assuming $s_3 \ne t_3$, the remaining graph has exactly two odd vertices,
namely $s_3$ and $t_3$. Hence, these vertices are in the
same connected component. However, once we no longer regard $s_3$ and $t_3$ as terminals and try to solve
the four terminal instance, we lose the Eulerianess property.
There are some efficient algorithms (e.g. \cite{SP-78,shi-80,tho-04,tho-09})
for the case of two pairs of terminals (without Eulerianess assumption)
but no linear time bound seems to be known.

The proof of \refth{feasibility_criterion} presented in \cite{sch-03}
is rather simple but non-constructive.
Our main result is as follows:
\begin{theorem}
\label{th:main}
    An Eulerian instance of the problem~\refeq{three_pairs} can be checked for feasibility in $O(m)$ time.
    If the check turns out positive, the desired path collection
    can be constructed in $O(m)$ time.
\end{theorem}

\section{The Algorithm}

\subsection{Preliminaries}
\label{subsec:prelim}

This subsection describes some basic techniques for working with edge-disjoint paths.
If the reader is familiar with network flow theory,
this subsection may be omitted.

Suppose we are given an undirected graph $G$ and a pair of distinct nodes $s$ (\emph{source})
and $t$ (\emph{sink}) from $VG$. An \emph{$s$--$t$ cut} is a subset $U \subseteq VG$
such that $s \in U$, $t \notin U$.

Edge-disjoint path collections can be described in flow-like terms as follows.
Let $\vec G$ denote the digraph obtained from $G$ by replacing each edge with
a pair of oppositely directed arcs. A subset $F \subseteq A\vec G$ is called \emph{balanced}
if $\abs{F \cap \deltain(v)} = \abs{F \cap \deltaout(v)}$ holds for each $v \in VG - \set{s,t}$.
Consider the \emph{value} of~$F$ defined as follows:
$$
	\val{F} := \abs{F \cap \deltaout(s)} - \abs{F \cap \deltain(s)}.
$$
Proofs of upcoming \reflm{decomp}, \reflm{augment}, and \reflm{min_cut} are quite standard and hence omitted
(see, e.g.~\cite{FF-62,CLRS-01}).

\begin{lemma}
\label{lm:decomp}
    Each balanced arc set decomposes into a collection arc-disjoint $s$--$t$ paths $\calP_{st}$,
    a collection of $t$--$s$ paths $\calP_{ts}$, and a collection of cycles $\calC$.
    Each such decomposition obeys $\abs{\calP_{st}} - \abs{\calP_{ts}} = \val{F}$.
    Also, such a decomposition can be carried out in $O(m)$ time.
\end{lemma}

Obviously, for each collection $\calP$ of edge-disjoint $s$--$t$ paths in $G$ there exists
a balanced arc set of value $\abs{\calP}$. Vice versa, each balanced arc set $F$ in $\vec{G}$
generates at least $\val{F}$ edge-disjoint $s$--$t$ paths in $G$.
Hence, finding a maximum cardinality collection of edge-disjoint $s$--$t$ paths in $G$
amounts to maximizing the value of a balanced arc set.

Given a balanced set $F$, consider the \emph{residual digraph}
$\vec G_F := (VG, (A\vec G - F) \cup F^{-1})$,
where $F^{-1} := \setst{a^{-1}}{a \in F}$ and $a^{-1}$ denotes the arc reverse to $a$.
\begin{lemma}
\label{lm:augment}
    Let $P$ be an arc-simple $s$--$t$ path in $\vec G_F$.
    Construct the arc set $F'$ as follows:
    (i)~take set $F$;
    (ii)~add all arcs $a \in AP$ such that $a^{-1} \notin F$;
    (iii)~remove all arcs $a \in F$ such that $a^{-1} \in AP$.
    Then, $F'$ is balanced and obeys $\val{F'} = \val{F} + 1$.
\end{lemma}
\begin{lemma}
\label{lm:min_cut}
    Suppose there is no $s$--$t$ path in $\vec G_F$. Then $F$ is of maximum value.
    Moreover, the set $U$ of nodes that are reachable from $s$ in~$\vec G_F$ obeys
    $d_G(U) = \val{F}$. Additionally, $U$ is an inclusion-wise minimum such set.
\end{lemma}

Hence, to find a collection of $r$ edge-disjoint $s$--$t$ paths one
needs to run a reachability algorithm in a digraph at most $r$ times. Totally, this
takes $O(rm)$ time and is known as the \emph{method of Ford and Fulkerson}~\cite{FF-62}.

\subsection{Checking for Feasibility}
\label{subsec:feasibility}

We start with a number of easy observations.
Firstly, there are some simpler versions of~\refeq{three_pairs}.
Suppose only one pair of terminals is given, i.e. $H = (VG, \set{s_1t_1})$.
Then the problem consists in checking if $s_1$ and $t_1$ are in the same
connected component of $G$. Note that if the instance $(G, H)$ is Eulerian
then it is always feasible since a connected component cannot contain
a single odd vertex. An $s_1$--$t_1$ path $P_1$ can be found in $O(m)$ time.

Next, consider the case of two pairs of terminals, i.e. $H = (VG, \set{s_1t_1, s_2t_2})$.
Connected components of $G$ not containing any of the terminals may be ignored.
Hence, one may assume that $G$ is connected since otherwise the problem
reduces to a pair of instances each having a single pair of terminals.

\begin{lemma}
\label{lm:two_pairs}
    Let $(G, H)$ be an Eulerian instance with two pairs of terminals.
    If $G$ is connected then $(G, H)$ is feasible.
    Also, the desired path collection $\set{P_1, P_2}$ can be found in $O(m)$ time.
\end{lemma}
\begin{proof}
    The argument is the same as in \refsec{intro}.
    Consider an arbitrary $s_1$--$t_1$ path $P_1$ and remove it from $G$.
    The resulting graph $G'$ may lose connectivity, however, $s_2$ and $t_2$
    are the only odd vertices in it (assuming $s_2 \ne t_2$).
    Hence, $s_2$ and $t_2$ are in the same
    connected component of $G'$, we can trace an $s_2$--$t_2$ path $P_2$ and, hence,
    solve the problem. The time complexity of this procedure is obviously $O(m)$.
\end{proof}

Now we explain how the feasibility of an Eulerian instance $(G, H)$
having three pairs of terminals can be checked in linear time.
Put $T := \set{s_1, s_2, s_3, t_1, t_2, t_3}$. There are exponentially many
subsets $U \subseteq VG$. For each subset $U$ consider its \emph{signature} $U^* := U \cap T$.
Fix an arbitrary signature $U^* \subseteq T$ and assume
w.l.o.g. that $\delta_H(U^*) = \set{s_1t_1, \ldots, s_kt_k}$.
Construct a new undirected graph $G(U^*)$ as follows:
add source~$s^*$,  sink~$t^*$, and $2k$ auxiliary edges $s^*s_1, \ldots, s^*s_k, t_1t^*, \ldots, t_kt^*$ to $G$.

Let $\nu(U^*)$ be the maximum cardinality of a collection of
edge-disjoint $s^*$--$t^*$ paths in $G(U^*)$.
We restate \refth{feasibility_criterion} as follows:
\begin{lemma}
\label{lm:feasibility_criterion_restated}
    An Eulerian instance $(G,H)$ with three pairs of terminals
    is feasible if and only if $\nu(U^*) \ge d_H(U^*)$ for each $U^* \subseteq T$,
\end{lemma}
\begin{proof}
    Necessity being obvious, we show sufficiency.
    Let $(G,H)$ be infeasible, then by \refth{feasibility_criterion}
    $d_G(U) < d_H(U)$ for some $U \subseteq VG$. Consider the corresponding signature $U^* := U \cap T$.
    One has $d_H(U) = d_H(U^*)$, hence there is a collection of $d_H(U)$ edge-disjoint
    $s^*$--$t^*$ paths in $G(U^*)$. Each of these paths crosses the cut $\delta_G(U)$ by a unique edge,
    hence $d_G(U) \ge d_H(U)$~--- a contradiction.
\end{proof}

By the above lemma, to check $(G,H)$ for feasibility one
has to validate the inequalities $\nu(U^*) \ge d_H(U^*)$ for all $U^* \subseteq T$.
For each fixed signature $U^*$ we consider graph $G(U^*)$,
start with an empty balanced arc set and perform up to three augmentations,
as explained in \refsubsec{prelim}. Therefore, the corresponding inequality is checked in $O(m)$ time.
The number of signatures is $O(1)$, which gives the linear
time for the whole feasibility check.

\medskip

We now present our first $O(m^2)$ time algorithm for finding the required path collection.
It will not be used in the sequel but gives some initial insight on the problem.
Consider an instance $(G, H)$ and let $s_1, t_1 \in T$ be a pair of terminals ($s_1t_1 \in EH$).
If $s_1 = t_1$ then the problem reduces to four terminals and is solvable in linear time,
as discussed earlier.
Otherwise, let $e$ be an edge incident to $s_1$, and put $s_1'$ to be the other endpoint of $e$.
We construct a new problem instance
$(G_e, H_e)$, where $G_e = G - e$, $H_e = H - s_1t_1 + s_1't_1$
(i.e. we remove edge $e$ and choose $s_1'$ instead of $s_1$ as a terminal).
Switching from $(G,H)$ to $(G_e,H_e)$ is called \emph{a local move}.
Local moves preserve Eulerianess.
If a local move generates a feasible instance then it is called \emph{feasible},
\emph{infeasible} otherwise.

If $(G_e,H_e)$ is feasible (say, $\calP_e$ is a solution)
then so is $(G,H)$ as we can append the edge $e$
to the $s_1'$--$t_1$ path in $\calP_e$ and obtain a solution $\calP$ to $(G,H)$.

Since $(G,H)$ is feasible, there must be a feasible local move (e.g. along the first
edge of an $s_1$--$t_1$ path in a solution).
We find this move by enumerating all edges $e$ incident to $s_1$ and checking
$(G_e,H_e)$ for feasibility. Once $e$ is found, we recurse to the instance $(G_e, H_e)$ having one
less edge. This way, a solution to the initial problem is constructed.

To estimate the time complexity, note that if a move along some edge $e$ is discovered to be
infeasible at some stage then it remains infeasible for the rest of the algorithm
(since the edge set of $G$ can only decrease). Hence, each edge in $G$ can be checked
for feasibility at most once. Each such check costs $O(m)$ time, thus yielding the total bound of $O(m^2)$.
This is already an improvement over the algorithm of Robertson and Seymour. However,
we can do much better.

\subsection{Reduction to a Critical Instance}

To solve an Eulerian instance of~\refeq{three_pairs} in linear time we start by constructing
an arbitrary node-simple $s_1$--$t_1$ path $P_1$.
Let $e_1, \ldots, e_k$ be the sequence of edges of $P_1$.
For each $i = 0, \ldots, k$ let $(G_i, H_i)$ be the instance obtained from the initial one $(G,H)$ by
a sequence of local moves along the edges $e_1, \ldots, e_i$.
In particular, $(G_0, H_0) = (G,H)$.

If $(G_k,H_k)$ is feasible (which can be checked in $O(m)$ time)
then the problem reduces to four terminals
and can be solved in linear time by~\reflm{two_pairs}.
Otherwise we find an index $j$ such that $(G_j,H_j)$ is a feasible
instance whereas $(G_{j+1},H_{j+1})$ is not feasible.

This is done by walking in the reverse direction along $P_1$ and considering
the sequence of instances $(G_k,H_k)$, \ldots, $(G_0, H_0)$.
Fix an arbitrary signature $U^*$ in $(G_k, H_k)$.
As we go back along $P_1$, terminal $s_1$ is moving. We apply these
moves to $U^*$ and construct a sequence of signatures $U_i^*$ in $(G_i, H_i)$.
($i = 1, \ldots, k$; in particular, $U_k^* = U^*$).
Let $\nu_i(U^*)$ be the maximum cardinality of an edge-disjoint collection of $s^*$--$t^*$
paths in $G_i(U^*_i)$.

Consider a consequent pair $G_{i+1}(U_{i+1}^*)$ and $G_i(U_i^*)$.
When $s_1$ is moved from node $v$ back to $v'$, edge $s^*v$
is removed and edges $s^*v'$ and $v'v$ are inserted.
Note, that this cannot decrease the maximum cardinality of an edge-disjoint $s^*$--$t^*$ paths collection
(if the dropped edge $s^*v$ was used by some path in a collection, then
we may replace it by a sequence of edges $s^*v'$ and $v'v$). Hence,
$$
    \nu_0(U_0^*) \ge \nu_1(U_1^*) \ge \ldots \ge \nu_k(U_k^*).
$$

Our goal is to find, for each choice of $U^*$, the largest
index $i$ (denote it by $j(U^*)$) such that
$$
    \nu_i(U_i^*) \ge d_H(U_i^*).
$$
Then, taking
$$
    j := \min_{U^* \subseteq T} j(U^*)
$$
we get the desired feasible instance $(G_j,H_j)$ such that $(G_{j+1},H_{j+1})$ is infeasible.

\medskip

To compute the values $\nu_i(U_i^*)$ consider the following dynamic problem.
Suppose one is given an undirected graph $\Gamma$ with distinguished
source $s$ and sink $t$, and also an integer $r \ge 1$.
We need the following operations:
\begin{quote}
    \textsc{Query}: Report $\min(r,c)$, where $c$ is the maximum cardinality
    of a collection of edge-disjoint $s$--$t$ paths in $\Gamma$.
\end{quote}
\begin{quote}
    \textsc{Move}$(v,v')$: Let $v$, $v'$ be a pair nodes in $VG$, $v \ne s$, $v' \ne s$, $sv \in E\Gamma$.
    Remove the edge $sv$ from $\Gamma$ and add edges $sv'$ and $v'v$ to $\Gamma$.
\end{quote}

\begin{lemma}
\label{lm:dyn_flows}
    There exists a data structure that can execute any sequence of \textsc{Query} and
    \textsc{Move} requests in $O(rm)$ time.
\end{lemma}
\begin{proof}
    We use a version of a folklore incremental reachability data structure.
    When graph $\Gamma$ is given to us, we start computing a balanced arc set $F$ in $\vec \Gamma$
    of maximum value $\val{F}$ but stop if $\val{F}$ becomes equal to $r$. This takes $O(rm)$ time.
    During the usage of the required data structure,
	the number of edge-disjoint $s$--$t$ paths (hence, $\val{F}$) cannot decrease
	(it can be shown using arguments similar to the described earlier).
    Therefore, if $\val{F} = r$ we may stop any maintenance and report the value of $r$ on each \textsc{Query}.

    Otherwise, as $\val{F}$ is maximum, there is no $s$--$t$
    path in $\vec \Gamma_F$ by \reflm{min_cut}.
    As long as no $r$ edge-disjoint $s$--$t$ paths in $\Gamma$ exist,
    the following objects are maintained:
    \begin{itemize}
        \item a balanced subset $F \subseteq A\Gamma$ of maximum value $\val{F}$ (which is less than $r$);
	\item an inclusionwise maximal directed tree $\calT$ rooted at $t$ and consisting of arcs
	from $\vec \Gamma_F$ (oriented towards to $t$).
    \end{itemize}

    In particular, $\calT$ covers exactly the set of nodes that can reach $t$
    by a path in $\vec \Gamma_F$. Hence, $s$ is not covered by $\calT$ (by~\reflm{augment}).

    Consider a \textsc{Move}$(v,v')$ request. We update $F$ as follows.
    If $sv \notin F$, then no change is necessary.
    Otherwise, we remove $sv$ from $F$ and also add arcs $sv'$ and $v'v$ to $F$.
    This way, $F$ remains balanced and $\val{F}$ is preserved.

    Next, we describe the maintenance of $\calT$.
    Adding an arbitrary edge $e$ to $\Gamma$ is simple.
    Recall that each edge in $\Gamma$ generates a pair of oppositely directed
    arcs in $\vec\Gamma$.
    Let $a = pq$ be one of these two arcs generated by $e$.
	Suppose $a$ is present in $\vec{\Gamma}_F$.
    If $a \in \deltain(V\calT)$ (i.e., $p$ is not reachable and $q$ is reachable)
    then add $a$ to $\calT$. Next, continue growing $\calT$ incrementally from $p$
    by executing a depth-first search and stopping at nodes already covered by $\calT$.
    This way, $\calT$ is extended to a maximum directed tree rooted at $t$.
    In other cases ($a \notin \deltain(V\calT)$) arc $a$ is ignored.

    Next consider deleting edge $sv$ from $G$. We argue that its removal
    cannot invalidate $\calT$, that is, $sv$ does not generate an arc from $\calT$.
    This is true since $t$ is not reachable from $s$ and, hence, arcs incident to $s$ may not appear in $\calT$.

    Note that a \emph{breakthrough} may occur during the above incremental procedure, i.e.
    node $t$ may become reachable from $s$ at some point.
    When this happens, we trace the corresponding $s$--$t$ path in $\calT$,
    augment $F$ according to \reflm{augment}, and recompute $\calT$ from scratch.
    Again, any further activity stops once $\val{F}$ reaches~$r$.

    To estimate the complexity, note that between breakthroughs
    we are actually dealing with a single suspendable depth-first traversal of $\vec \Gamma_F$.
    Each such traversal costs $O(m)$ time and there are at most $r$ breakthroughs.
    Hence, the total bound of $O(rm)$ follows.
\end{proof}

We apply the above data structure to graph $G_k(U_k^*)$ for $r = d_H(U^*)$
and make the moves in the reverse order, as explained earlier.
Once \textsc{Query} reports the existence
of $r$ edge-disjoint $s^*$--$t^*$ paths in $G_i(U_i^*)$, put $j(U^*) := i$ and proceed
to the next signature. This way, each value $j(U^*)$ can be computed in $O(m)$ time.
There are $O(1)$ signatures and $r = O(1)$, hence computing $j$ takes linear time as well.

\subsection{Dealing with a Critical Instance}

The final stage deals with problem instance $(G_j, H_j)$.
For brevity, we reset $G := G_j$, $H := H_j$ and also denote $G' := G_{j+1}$, $H' := H_{j+1}$.
Consider the connected components of $G$.
Components not containing any terminals may be dropped.
If $G$ contains at least two components with terminals,
the problem reduces to simpler cases described in \refsubsec{feasibility}.
Hence, one can assume that $G$ is connected.
We prove that $(G,H)$ is, in a sense, \emph{critical}, that is,
it admits a cut of a very special structure.

\begin{figure}[t!]
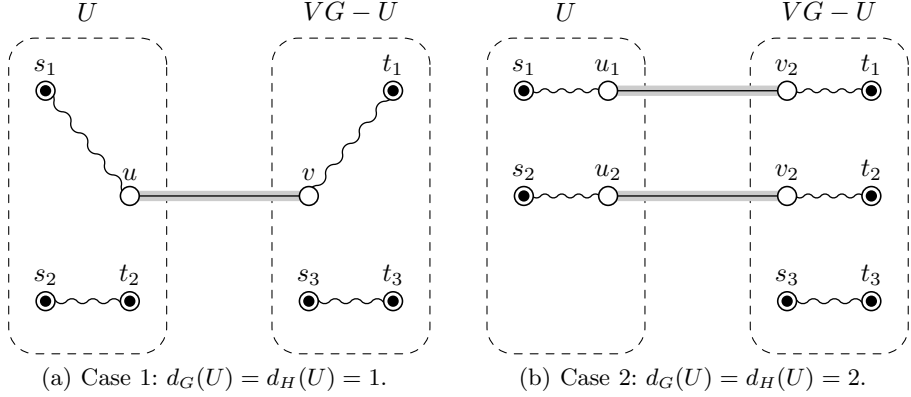

    \centering
    \subfigure[Case~1: $d_G(U) = d_H(U) = 1$.]{
        \includegraphics{pics/critical.1}%
    }
    \hspace{0.5cm}%
    \subfigure[Case~2: $d_G(U) = d_H(U) = 2$.]{
        \includegraphics{pics/critical.2}%
    }
    \caption{
        A critical instance $(G,H)$.
        Terminals are marked with dots and renumbered.
        Wavy lines indicate parts of paths in the desired collection $\calP$.
    }
\label{fig:critical}
\end{figure}

\begin{lemma}
\label{lm:critical_cut}
    There exists a subset $U \subseteq VG$ such that $d_G(U) = d_H(U) = 2$, $G[U]$ is connected and $\abs{U \cap T} = 2$
    (see~\reffig{critical}(b)).
\end{lemma}
\begin{proof}
    The following is true:
    \begin{numitem}
    \label{eq:step_cond}
        For problem instances $(G,H)$ and $(G',H')$
        \begin{itemize}
            \item[(i)] $(G,H)$ is feasible,
            \item[(ii)] $(G',H')$ is obtained from $(G,H)$ by a single local move,
            \item[(iii)] $(G',H')$ is infeasible,
            \item[(iv)] let $s \in VG'$ be the new location of the moved terminal, $st \in EH'$,
            then $s$ and $t$ are in the same connected component of $G'$.
        \end{itemize}
    \end{numitem}

    Properties (i)--(iii) are ensured by the choice of $j$. Property (iv) holds since there exists
    a remaining (untraversed) part of the initial $s_1$--$t_1$ path in the original graph $G$.

    We apply a number of contractions to $(G,H)$ and $(G',H')$ that preserve condition~\refeq{step_cond}.
    Suppose the following:
    \begin{numitem}
    \label{eq:bad_bridge}
        there is a subset $U \subseteq VG$ such that $d_G(U) = d_H(U) = 1$
        and $\abs{U \cap T} = 1$.
    \end{numitem}
    In other words, there is a \emph{bridge} $e = uv \in EG$, $u \in U$, $v \in VG - U$
    (an edge whose removal increases the number of connected components)
    that separates $G$ into parts $G[U]$ and $G[VG - U]$ and the former part
    contains a single terminal, say $s$.

    We argue that the local move, which produced $(G',H')$,
    was carried out in the subgraph $G[VG - U]$
    (but not in $G[U]$ or along the edge $e$).

    Firstly, the move could not have been applied to $s$.
    Suppose the contrary.
    Terminal $s$ is connected to node $v$ by some path in $G[U \cup \set{v}]$
    and this property remains true even if apply a local move to $s$.
    (Nodes $v$ and $s$ are the only odd vertices in $G[U \cup \set{v}]$,
    hence, these nodes cannot fall into distinct connected components after the move.)
    Therefore, $(G,H)$ and $(G',H')$ are simultaneously feasible or infeasible.

    Next, suppose that $v$ is a terminal and the move is carried out along the bridge $e$.
    Then, $vs \notin EH'$ (otherwise, $(G',H')$ remains feasible).
    Therefore, $vw \in EH'$ for some $w \in VG - U$.
    Then $v$ and $w$ belong to different connected components of $G'$ after the move,
    which is impossible by \refeq{step_cond}(iv).

    Contract the set $U \cup \set{v}$ in instances $(G, H)$ and $(G', H')$ thus producing
    instances $(\bar G, \bar H)$ and $(\bar G', \bar H')$, respectively.
    The above contraction preserves feasibility, hence $(\bar G, \bar H)$ is feasible
    and $(\bar G', \bar H')$ is infeasible. Moreover, the latter instance is obtained from the former one
    by a local move. Property~\refeq{step_cond} is preserved.

    \medskip

    We proceed with these contractions until no subset $U$ obeying \refeq{bad_bridge} remains.
    Next, since $(G',H')$ is infeasible by \refth{feasibility_criterion} there exists
    a subset $U \subseteq VG$ such that $d_{G'}(U) < d_{H'}(U)$.
    Eulerianess of $G' + H'$ implies that each cut in $G' + H'$ is crossed
    by an even number of edges, hence $d_{G'}(U) \equiv d_{H'}(U) \pmod{2}$.
    Therefore,
    \begin{equation}
    \label{eq:bad_cut}
        d_{G'}(U) \le d_{H'}(U) - 2.
    \end{equation}
    At the same time, $(G,H)$ is feasible and hence
    \begin{equation}
    \label{eq:good_cut}
        d_G(U) \ge d_H(U).
    \end{equation}
    Graph $G'$ is obtained from $G$ by removing a single edge.
    Similarly, $H'$ is obtained from $H$ by one edge removal and one edge insertion.
    Hence, $d_G(U)$ and $d_H(U)$ differ from $d_{G'}(U)$ and $d_{H'}(U)$ (respectively)
    by at most~1. Combining this with \refeq{bad_cut} and \refeq{good_cut}, one has
    $$
        d_{G'}(U) + 1 = d_G(U) = d_H(U) = d_{H'}(U) - 1.
    $$
    So $d_H(U) \in \set{1,2}$.

    Suppose $d_H(U) = 1$.
    Subgraphs $G[U]$ and $G[VG - U]$ are connected (since otherwise $G$
    is not connected). Also, $\abs{U \cap T} = 3$ (otherwise, $\abs{U \cap T} = 1$
    or $\abs{(VG - U) \cap T} = 1$ and \refeq{bad_bridge} still holds).
    Therefore, Case~1 from \reffig{critical}(a) applies
    (note that terminals $s_i$ and $t_i$ depicted there are appropriately renumbered).
    Let us explain, why this case is impossible.
    Graph $G'$ is obtained from $G$ by removing edge $uv$.
    Let, as in \refeq{step_cond}(iv), $s$ denote the terminal in $(G,H)$
    that is being moved and $t$ denote its ``mate'' terminal (i.e. $st \in EH$).
    We can assume by symmetry that $u = s$.
    Hence, $v$ is the new location of $s$ in $(G',H')$.
    By \refeq{step_cond}(iv), $v$ and $t$ are in the same connected component of $G'$. 
    The latter is only possible if $s = u = s_1$ and $t = t_1$.
    But then feasibility of $(G,H)$ implies that of $(G',H')$.

    Finally, let $d_H(U) = 2$. Replacing $U$ by $VG - U$, if necessary,
    we may assume that $\abs{U \cap T} = 2$, see \reffig{critical}(b).
    It remains to prove that $G[U]$ is connected.
    Let us assume the contrary. Then, $U = U_1 \cup U_2$, $U_1 \cap U_2 = \emptyset$,
    $d_G(U_1) = d_H(U_1) = 1$, $d_G(U_2) = d_H(U_2) = 1$ (due to feasibility of $(G,H)$ and connectivity of $G$).
    Therefore, \refeq{bad_bridge} still holds (both for $U := U_1$ and $U := U_2$)~--- a contradiction.

    Once set $U$ is found, we undo the contractions described in the beginning and
    obtain a set $U$ for the original instance $(G,H)$. Clearly, these uncontractions preserve
    the required properties of $U$.
\end{proof}

\begin{lemma}
\label{lm:building_u}
    Set $U$ figured in \reflm{critical_cut} can be constructed in $O(m)$ time.
\end{lemma}
\begin{proof}
    We enumerate pairs of terminals $p, q \in T$ that might qualify for $U^* := U \cap T = \set{p,q}$.
    Take all such pairs $U^* = \set{p,q}$ except those forming an edge in~$H$ ($pq \in EH$).
    Contract $U^*$ and $T - U^*$ into $s^*$ and $t^*$, respectively, in the graphs $G$ and $H$.
    The resulting graphs are denoted by $G^*$ and $H^*$.
    If a subset obeying \reflm{critical_cut} and having the signature $U^*$ exists then there must be an $s^*$--$t^*$
    cut $U$ in $G^*$ such that $d_{G^*}(U) = 2$.

    We try to construct $U$ by applying three iterations of
    the max-flow algorithm of Ford and Fulkerson, see \refsubsec{prelim}.
    If the third iteration succeeds, i.e. three edge-disjoint $s^*$--$t^*$ paths
    are found, then no desired cut $U$ having signature $U^*$ exists;
    we continue with the next choice of $U^*$.
    Otherwise, a subset $U \subseteq VG^*$ obeying $d_{G^*}(U) \le 2$ is constructed.
    Case $d_{G^*}(U) < 2 = d_{H^*}(U)$ is not possible due to feasibility of $(G,H)$.

    Set $U$ is constructed for graph $G^*$ but may also be regarded
    as a subset of $VG$. We keep notation $U$ when referring to this subset.

    Connectivity of $G[U]$ is the only remaining property we need to ensure.
    This is achieved by selecting an inclusion-wise maximal set $U$ among minimum-capacity cuts
    that separate $\set{p,q}$ and $T - \set{p,q}$.
    Such maximality can achieved in a standard way, i.e. by traversing the residual network
    in backward direction from the sink~$t^*$, see \reflm{min_cut}.

    To see that $G[U]$ is connected suppose the contrary.
    Then, as in the proof of \reflm{critical_cut},
    let $U_1$ and $U_2$ be the node sets of the connected components of $G[U]$.
    Edges in $\delta_G(U) = \set{e_1, e_2}$ are bridges connecting $G[U_i]$ to the remaining
    part of graph $G$ (for $i = 1, 2$). Also, $\abs{U_1 \cap T} = \abs{U_2 \cap T} = 1$
    (recall that $G$ is connected).
    Denote $U_1 \cap T = \set{q_1}$ and $U_2 \cap T = \set{q_2}$.
    Terminals $q_1$ and $q_2$ are not connected in $G[U]$.
    Since set~$U$ is inclusion-wise maximal, any subset $U'$ satisfying \reflm{critical_cut} also obeys $U' \subseteq U$.
    But then $q_1$ and $q_2$ are also disconnected in $G[U']$, which is a contradiction.
    Therefore, no valid subset $U$ of signature $U^*$ obeying \reflm{critical_cut} exists.

    In the algorithm, we check $G[U]$ for connectivity in $O(m)$ time.
    If the graph is not connected, then we proceed with the next signature $U^*$.
\end{proof}

Now everything is in place to complete the proof of \refth{main}.
By \reflm{building_u}, finding set~$U$ takes $O(m)$ time.
It remains to construct a solution to $(G,H)$.
Put $\delta_G(U) = \set{e_1, e_2}$, $e_i = u_iv_i$, $u_i \in U$, $v_i \in VG - U$,
$i = 1, 2$. Again, after renaming some terminals we may assume that
$s_1, s_2 \in U$, $t_1, t_2, s_3, t_3 \in VG - U$.
Augment $G$ by adding two new nodes $s^*$ and $t^*$ and \emph{auxiliary} edges
$s^*u_1$, $s^*u_2$, $t_1t^*$, and $t_2t^*$.
Due to feasibility of $(G,H)$, there exists (and can be constructed in $O(m)$ time)
a collection of two edge-disjoint $s^*$--$t^*$ paths.
After removing auxiliary edges we either obtain
a $u_1$--$t_1$ path and a $u_2$--$t_2$ path (Case A)
or a $u_1$--$t_2$ path and a $u_2$--$t_1$ path (Case B).
To extend these paths to an $s_1$--$t_1$ path and an $s_2$--$t_2$ path
we consider a four terminal instance in the subgraph $G[U]$.
The demand graph is $(U, \set{s_1u_1, s_2u_2})$ in Case A and $(U, \set{s_1u_2, s_2u_1})$ in Case B.
As $G[U]$ is connected, the latter instance is feasible by \reflm{two_pairs}.
Therefore, we obtain edge-disjoint $s_1$--$t_1$ and $s_2$--$t_2$ paths $P_1$ and $P_2$, respectively.
As explained earlier in \refsec{intro}, the remaining path $P_3$ always exists
and can be found in $O(m)$ time.
Therefore, the proof of \refth{main} is complete.

\nocite{*}
\bibliographystyle{alpha}
\bibliography{main}

\end{document}